\documentclass[leqno,12pt]{article}

\usepackage{amssymb}
\usepackage{euscript}
\usepackage[dvips]{graphicx}
\usepackage{flafter}
\usepackage{pstricks}
\usepackage[latin1]{inputenc}
\usepackage{amsmath}
\usepackage{amsfonts}
\usepackage{pstricks-add}

\usepackage{mathrsfs} 

\evensidemargin\oddsidemargin

\begin{document}

\baselineskip=18pt \setcounter{page}{1}

\renewcommand{\theequation}{\thesection.\arabic{equation}}
\newtheorem{theorem}{Theorem}[section]
\newtheorem{lemma}[theorem]{Lemma}
\newtheorem{proposition}[theorem]{Proposition}
\newtheorem{corollary}[theorem]{Corollary}
\newtheorem{remark}[theorem]{Remark}
\newtheorem{fact}[theorem]{Fact}
\newtheorem{problem}[theorem]{Problem}
\newtheorem{example}[theorem]{Example}
\newtheorem{question}[theorem]{Question}
\newtheorem{conjecture}[theorem]{Conjecture}

\newcommand{\eqnsection}{
\renewcommand{\theequation}{\thesection.\arabic{equation}}
    \makeatletter
    \csname  @addtoreset\endcsname{equation}{section}
    \makeatother}
\eqnsection

\def\r{{\mathbb R}}
\def\e{{\mathbb E}}
\def\p{{\mathbb P}}
\def\P{{\bf P}}
\def\E{{\bf E}}
\def\Q{{\bf Q}}
\def\z{{\mathbb Z}}
\def\N{{\mathbb N}}
\def\T{{\mathbb T}}
\def\G{{\mathbb G}}
\def\L{{\mathbb L}}

\def\deg{\chi}

\def\ee{\mathrm{e}}
\def\d{\, \mathrm{d}}
\def\S{\mathscr{S}}




\vglue50pt

\centerline{\bf Weak convergence for the minimal position in a branching random walk:}

\medskip
\centerline{\bf a simple proof}

\bigskip
\bigskip

\centerline{by}

\medskip

\centerline{Elie A\"\i d\'ekon $\;$and$\;$ Zhan Shi}

\medskip

\centerline{\it Technische Universiteit Eindhoven $\;$\&$\;$ Universit\'e Paris VI}

\bigskip

\centerline{\it Dedicated to Professors Endre Cs\'aki and P\'al R\'ev\'esz}

\centerline{\it on the occasion of their 75th birthdays}

\bigskip
\bigskip

{\leftskip=2truecm \rightskip=2truecm \baselineskip=15pt \small

\noindent{\slshape\bfseries Summary.} Consider the boundary case in a one-dimensional super-critical branching random walk. It is known that upon the survival of the system, the minimal position after $n$ steps behaves in probability like ${3\over 2} \log n$ when $n\to \infty$. We give a simple and self-contained proof of this result, based exclusively on elementary properties of sums of i.i.d.\ real-valued random variables.

\bigskip


\bigskip


} 

\bigskip
\bigskip

\section{Introduction}
   \label{s:intro}

$\phantom{aob}$Consider a (discrete-time, one-dimensional) branching random walk. It starts with an initial ancestor particle located at the origin. At time $1$, the particle dies, producing a certain number of new particles; these new particles are positioned according to the law of a given finite point process. At time $2$, these particles die, each giving birth to new particles that are positioned (with respect to the birth place) according to the law of the same point process. And the system goes on indefinitely, as long as there are particles that are alive. We assume that, each particle produces new particles independently of other particles in the same generation, and of everything up to that generation.

The number of particles in each generation obviously forms a Galton--Watson process, which will always be assumed to be super-critical.

Let $(V(x), \; |x|=n)$ be the positions of the particles at the $n$-th generation. The process $(V(x))$ indexed by a Galton--Watson tree is called a branching random walk. We do not assume the random variables $V(x)$, $|x|=1$, to be independent, nor necessarily identically distributed, though it is often assumed in the literature (for example, in \cite{R94}).

We are interested in $\min_{|x|=n} V(x)$, the minimal position of the branching random walk after $n$ steps. Under a mild integrability assumption, we have (Kingman~\cite{kingman}, Hammersley~\cite{hammersley}, Biggins~\cite{biggins}), on the set of non-extinction,
\begin{equation}
    {1\over n} \min_{|x|=n} V(x)
    \to
    \gamma, \qquad\hbox{\rm a.s.},
    \label{LGN}
\end{equation}

\noindent where $\gamma\in \r$ is a known constant. 

The rate of convergence in (\ref{LGN}) has recently been studied, independently, by Hu and Shi~\cite{yzpolymer}, and Addario-Berry and Reed~\cite{addario-berry-reed}. To state the result, we assume the following condition:
\begin{equation}
    \E \Big( \sum_{|x|=1} \ee^{-V(x)} \Big) =1,
    \qquad
    \E \Big( \sum_{|x|=1} V(x) \ee^{-V(x)} \Big) =0.
    \label{cond-hab}
\end{equation}

\noindent This is referred to in the literature as the boundary case; see for example Biggins and Kyprianou~\cite{biggins-kyprianou05}. Under (\ref{cond-hab}), we have $\gamma=0$ in (\ref{LGN}).

For discussions on the nature of the assumption (\ref{cond-hab}), see Jaffuel~\cite{jaffuel}. Loosely speaking, letting $\underline{m}$ denote the essential infimum of $\min_{|x|=1}V(x)$, then under some mild integrability conditions, a branching random walk can always be made to satisfy (\ref{cond-hab}) after a suitable change of scale, if either $\underline{m}=-\infty$, or $\underline{m}>-\infty$ and $\E[\sum_{|x|=1} {\bf 1}_{ \{ x= \underline{m}\} }]<1$.

In addition of (\ref{cond-hab}), we assume the following integrability condition: there exists $\delta>0$ such that
\begin{equation}
    \E \Big( \sum_{|x|=1} \ee^{\delta V(x)} \Big) <\infty,
    \quad
    \E \Big( \sum_{|x|=1} V(x) \ee^{-(1+\delta)V(x)} \Big) <\infty,
    \quad
    \E \Big[ \Big( \sum_{|x|=1} 1\Big)^{1+\delta} \, \Big] <\infty.
    \label{delta}
\end{equation}

\medskip

\begin{theorem}
 \label{t:main}
 {\bf (\cite{yzpolymer}, \cite{addario-berry-reed})}
 Assume (\ref{cond-hab}) and (\ref{delta}). On the set of non-extinction, we have
 \begin{equation}
     {1\over \log n} \min_{|x|=n} V(x)
     \to {3\over 2},
     \qquad\hbox{in probability.}
     \label{main}
 \end{equation}

\end{theorem}

\medskip

The proofs of Theorem \ref{t:main} presented in \cite{yzpolymer} and \cite{addario-berry-reed} are totally different, but both of them are rather technical. [In \cite{addario-berry-reed}, it is furthermore assumed that the random variables $V(x)$, $|x|=1$, are i.i.d.\ and that $\sum_{|x|=1}1$ is a.s.\ bounded.]

The sole goal of this note is to give a simple and self-contained proof of Theorem \ref{t:main}, using only elementary properties of sums of i.i.d.\ real-valued random variables. We list in the Appendix these elementary properties of sums of i.i.d.\ real-valued random variables.

Theorem \ref{t:main} is equivalent to saying that for any $\alpha<{3\over 2}$ and $\beta>{3\over 2}$,
\begin{eqnarray}
    \P\Big\{ \min_{|x|=n} V(x) \le \alpha\log n \, \Big| \, \hbox{non-extinction}\Big\} 
 &\to& 0,
    \label{lb}
    \\
    \P\Big\{ \min_{|x|=n} V(x) \ge \beta\log n \, \Big| \, \hbox{non-extinction}\Big\} 
 &\to& 0, 
    \qquad n\to \infty.
    \label{ub}
\end{eqnarray}

\noindent We prove (\ref{lb}) and (\ref{ub}) in Sections \ref{s:lb} and \ref{s:ub}, respectively. Our proof of (\ref{ub}) is presented under an additional assumption: $\E\{ [\sum_{|x|=1}1 ]^2\} <\infty$ and $\E\{ [\sum_{|x|=1}\ee^{-(1+\delta)V(x)}]^2\} <\infty$ for some $\delta>0$. It is possible to adapt the proof without using the additional assumption, by means of a truncation argument (see, for example, Lemma 4.5 of Gantert et al.~\cite{nyzsurvival}) and at the cost of an extra page, but we think it is more interesting to keep the proof as simple as possible.

Throughout the paper, we use $a_n \sim b_n$ ($n\to \infty$) to denote $\lim_{n\to \infty} \, {a_n\over b_n} =1$; the letter $c$ with subscript denotes a finite and positive constant.

{\bf Remark.} In order to make our proof truly self-contained, we reprove all known results for branching random walks that are needed in the paper.

\section{Proof of (\ref{lb})}
\label{s:lb}

$\phantom{aob}$For any vertex $x$, let $[\![ \varnothing, \, x]\!]$ be the unique shortest path relating $x$ to the root $\varnothing$, and $x_i$ (for $0\le i\le |x|$) the vertex on $[\![ \varnothing, \, x]\!]$ such that $|x_i|=i$.

We assume (\ref{cond-hab}). Let $S_1$, $S_2-S_1$, $S_3-S_2$, $\cdots$ be i.i.d.\ random variables such that
\begin{equation}
    \P\{ S_1 \le u\}
    =
    \E\Big\{ \sum_{|x|=1} \ee^{-V(x)} \,
    {\bf 1}_{ \{ V(x) \le u\} } \Big\} ,
    \qquad \forall u\in \r.
    \label{S1}
\end{equation}

\noindent In particular, $\E(S_1)=0$. We write $\underline{S}_n := \min_{1\le i\le n} S_i$ for $n\ge 1$.

We claim that for any $n\ge 1$ and any measurable function $g: \r^n \to [0, \, \infty)$,
\begin{equation}
    \E\Big\{ \sum_{|x|=n} g(V(x_1), \cdots, V(x_n)) \Big\}
    =
    \E \Big\{ \ee^{S_n} g(S_1, \cdots, S_n)\Big\} .
    \label{many-to-one}
\end{equation}

\noindent This is easily checked\footnote{See the last Remark in the Introduction. There is a deep explanation to the presence of the new random walk $(S_i)$ using the size-biased branching random walk of Lyons, Pemantle and Peres~\cite{lyons-pemantle-peres} and Lyons~\cite{lyons}. This idea has been used by many authors in various forms, going back at least to Kahane and Peyri\`ere~\cite{kahane-peyriere}.} by induction (on $n$): For $n=1$, (\ref{many-to-one}) is nothing else but (\ref{S1}). Assume that (\ref{many-to-one}) is proved for $n$. Then by conditioning on the first generation of the branching random walk and using the induction hypothesis, we obtain the claimed equality for $n+1$.

\bigskip

\noindent {\it Proof of (\ref{lb}).} Let $K>0$ and $0<a< {3\over 2}$. Let
$$
Z_n := \sum_{|x|=n} {\bf 1}_{ \{ V(x) \le a \log n, \; V(x_i) \ge -K, \; \forall 1\le i\le n\} }.
$$

\noindent By (\ref{many-to-one}), we have,
$$
\E(Z_n) = \E \Big\{ \ee^{S_n} \, {\bf 1}_{ \{ S_n \le a \log n, \; S_i \ge -K, \; \forall 1\le i\le n\} } \Big\} \le n^a \, \P\Big\{ S_n \le a \log n, \; \underline{S}_n \ge -K \Big\} .
$$

\noindent For $n$ such that $a\log n \ge 1$, we have $\P \{ S_n \le a \log n, \; \underline{S}_n \ge -K\} \le c_1 \, {(\log n)^2\over n^{3/2}}$ (Lemma \ref{l:3}). Since $a<{3\over 2}$, it follows that $\lim_{n\to \infty} \E(Z_n) =0$.

Under the assumption $\E[\sum_{|x|=1} \ee^{-V(x)}]=1$, $(\sum_{|x|=n} \ee^{-V(x)}, \; n\ge 0)$ is a (non-negative) martingale with respect to its natural filtration; so it converges almost surely. In particular, $\sup_n \sum_{|x|=n} \ee^{-V(x)}<\infty$ a.s. A fortiori, $\inf_n \min_{|x|=n} V(x) > -\infty$ a.s.\footnote{See the last Remark in the Introduction. In fact, assumption (\ref{cond-hab}) ensures that $\min_{|x|=n} V(x) \to \infty$ a.s.\ on the set of non-extinction; see Lyons~\cite{lyons}.}

Since we have already proved that $Z_n \to 0$ in probability, this implies (\ref{lb}).\hfill$\Box$

\section{Proof of (\ref{ub})}
\label{s:ub}

$\phantom{aob}$The proof of (\ref{ub}) also relies on the study of the associated random walk $(S_i)$. It is technically slightly more involved, because we use a second-moment argument this time. We start with the observation\footnote{See the last Remark in the Introduction. It is, obviously, an immediate consequence of (\ref{LGN}).} that, under assumption (\ref{delta}), there exists a constant $c_2<\infty$ such that, on the set of non-extinction,
\begin{equation}
    \limsup_{n\to \infty} \, {1\over n} \max_{|x|=n} V(x)
    \le
    c_2, \qquad\hbox{\rm a.s.}
    \label{LGN2}
\end{equation}

\noindent Indeed, by (\ref{delta}), there exists $\delta>0$ such that $c_3:= \E[ \sum_{|x|=1} \ee^{\delta V(x)}]<\infty$. Let $a>{\log c_3\over \delta}$, and let $Y_n := \sum_{|x|=n} {\bf 1}_{ \{ V(x) \ge an\} }$. Clearly, $Y_n \le \ee^{-an \delta} \sum_{|x|=n} \ee^{\delta V(x)}$. By Chebyshev's inequality, $\P(Y_n > 0) \le \E(Y_n) \le \ee^{-an \delta} \E[\sum_{|x|=n} \ee^{\delta V(x)}] = \ee^{-an \delta} \, (c_3)^n$, which is summable in $n$ because $c_3 < \ee^{a \delta}$. By the Borel--Cantelli lemma, $\lim_{n\to \infty} Y_n =0$ a.s.\ on the set of non-extinction, yielding (\ref{LGN2}).

\bigskip

\noindent {\it Proof of (\ref{ub}).} Let $C>0$ be the constant in Lemma \ref{l:4}. Let
$$
I_k = I_k(n) :=
\begin{cases}
[0, \, \infty) &\hbox{ if }1\le k\le n,
\\
[{3\over 2}\log n, \, \infty) &\hbox{ if }n<k<2n,
\\
[{3\over 2}\log n, \; {3\over 2}\log n + 2C] &\hbox{ if }k=2n.
\end{cases}
$$

\noindent Consider the random variable
\begin{eqnarray*}
    Z_n
 &:=& \sum_{|x|=2n}
    {\bf 1}_{ \{ V(x) \in [{3\over 2}\log n, \, {3\over 2}\log n +2C],
    \; \min_{1\le i\le n} V(x_i) \ge 0, \;
    \min_{n<j\le 2n} V(x_j) \ge {3\over 2}\log n\} }
    \\
 &=& \sum_{|x|=2n}
    {\bf 1}_{ \{ V(x_k) \in I_k, \; \forall 1\le k\le 2n\} }.
\end{eqnarray*}

\noindent  Applying (\ref{many-to-one}) gives that
$$
\E(Z_n)
=
\E \Big\{ \ee^{S_{2n}} \,
    {\bf 1}_{ \{ S_k\in I_k, \; \forall 1\le k\le 2n\} } \Big\}
\ge n^{3/2} \,
    \P\Big\{ S_k \in I_k, \; \forall 1\le k\le 2n\Big\} .
$$

\noindent By Lemma \ref{l:4}, this yields $\E(Z_n) \ge c_4$ for some constant $c_4>0$.

We now estimate the second moment. By definition,
\begin{eqnarray*}
    \E(Z_n^2)
 &=& \E \Big\{
    \sum_{|x|=2n}\sum_{|y|=2n}
    {\bf 1}_{ \{ V(x_k), \, V(y_k) \in I_k, \;
    \forall 1\le k\le 2n\} } \Big\}
    \\
 &=& \E(Z_n) +
    \E \Big\{ \sum_{j=0}^{2n-1} \sum_{|z|=j}
    {\bf 1}_{ \{ V(z_i) \in I_i, \; \forall 1\le i\le j\} }
    \sum_{(x_{j+1}, \, y_{j+1})} \sum_{(x, \, y)}
    {\bf 1}_{ \{ V(x_k), \, V(y_k) \in I_k, \; \forall j< k\le 2n\} }
    \Big\}
    \\
 &=:& \E(Z_n) + \Lambda_n,
\end{eqnarray*}

\noindent where, the double sum $\sum_{(x_{j+1}, \, y_{j+1})}$ is over pairs $(x_{j+1},\, y_{j+1})$ of distinct children of $z$, whereas $\sum_{(x,\, y)}$ is over pairs $(x,\, y)$ with $|x|=|y|=2n$ such that\footnote{By $x\ge y$, we mean either $x=y$, or $y$ is an ancestor of $x$.} $x\ge x_{j+1}$ and $y\ge y_{j+1}$.

Applying the Markov property at generation $j+1$ gives that
\begin{eqnarray}
    \Lambda_n
 &=&\E \Big\{
    \sum_{j=0}^{2n-1} \sum_{|z|=j}
    {\bf 1}_{ \{ V(z_i)\in I_i, \; \forall 1\le i\le j\} } \times
    \nonumber
    \\
 && \qquad\qquad
    \times \sum_{(x_{j+1}, \, y_{j+1})}
    {\bf 1}_{ \{ V(x_{j+1}), \,
    V(y_{j+1}) \in I_{j+1}\} }
    \, f_{j,n}(V(x_{j+1}))
    f_{j,n}(V(y_{j+1}))  \Big\},
    \label{E(Zn2)<}
\end{eqnarray}

\noindent where, for any $u\in \r$,
$$
f_{j,n}(u)
:= \E \Big\{ \sum_{|x|=2n-j-1}
    {\bf 1}_{ \{ u+V(x_\ell) \in I_{\ell + j+1}, \;
    \forall 1\le \ell\le 2n-j-1\} } \Big\} .
$$

\noindent Applying (\ref{many-to-one}), we get:
\begin{eqnarray*}
    f_{j,n}(u)
 &=& \E \Big\{ \ee^{S_{2n-j-1}}
    {\bf 1}_{ \{ u+S_\ell \in I_{\ell + j+1}, \;
    \forall 1\le \ell\le 2n-j-1\} } \Big\}
    \label{app-many-to-one}
    \\
 &\le& n^{3/2} \ee^{2C-u}\,
    \P\Big\{ \underline{S}_{2n-j-1} \ge -u, \;
    {3\over 2} \log n - u \le S_{2n-j-1} \le
    {3\over 2} \log n - u +2C \Big\} ,
\end{eqnarray*}

\noindent where $\underline{S}_k := \min_{1\le i\le k} S_i$ as before. By (\ref{A1}),
$$
f_{j,n}(u) \le n^{3/2} \ee^{2C-u}\, \P\Big\{ {3\over 2} \log n - u \le S_{2n-j-1} \le {3\over 2} \log n - u +2C \Big\} \le c_5 \, {n^{3/2}\ee^{-u}\over (2n-j)^{1/2}} .
$$

\noindent This inequality turns out to be too rough when $j\le n$, so we do differently in the latter situation. Since $2n-j-1 \ge n-1$ this time, we apply Lemma \ref{l:3}: for $u\in I_{j+1}$,
\begin{eqnarray*}
    f_{j,n}(u)
 &\le& n^{3/2} \ee^{2C-u}\,
    \P\Big\{ \underline{S}_{2n-j-1} \ge -u, \;
    S_{2n-j-1} \le {3\over 2} \log n - u +2C \Big\}
    \\
 &\le& c_6 \, n^{3/2} \ee^{-u}
    {(u+1) (u + \log n)^2\over (2n-j-1)^{3/2}}
    \\
 &\le& c_6 \, \ee^{-u} (u+1) (u + \log n)^2
    \\
 &\le& c_7 \, (\log n)^2 \ee^{-u} (u+1)^3.
\end{eqnarray*}

Let us go back to (\ref{E(Zn2)<}), to see that
\begin{eqnarray*}
    \Lambda_n
 &\le&c_7^2 (\log n)^4\, 
    \E \Big\{ \sum_{j=0}^n \sum_{|z|=j}
    {\bf 1}_{ \{ V(z_i) \in I_i, \; \forall 1\le i\le j\} } \times
    \\
 &&\qquad \times 
    \sum_{(x_{j+1}, \, y_{j+1})} 
    \ee^{-V(x_{j+1})-V(y_{j+1})} 
    |V(x_{j+1})+1|^3\, |V(y_{j+1})+1|^3 \Big\}
    \\
 &&+ c_5^2 n^3\, 
    \E \Big\{ \sum_{j=n+1}^{2n-1} {1\over 2n-j}
    \sum_{|z|=j}
    {\bf 1}_{ \{ V(z_i) \in I_i, \; \forall 1\le i\le j\} }
    \sum_{(x_{j+1}, \, y_{j+1})}
    \ee^{-V(x_{j+1})-V(y_{j+1})} \Big\}.
\end{eqnarray*}

\noindent We observe that $\sum_{(x_{j+1}, \, y_{j+1})} \ee^{-V(x_{j+1})-V(y_{j+1})} |V(x_{j+1})+1|^3\, |V(y_{j+1})+1|^3$ is bounded by $[\sum_{x_{j+1}} \ee^{-V(x_{j+1})} |V(x_{j+1})+1|^3]^2$, whereas $\sum_{(x_{j+1}, \, y_{j+1})} \ee^{-V(x_{j+1})-V(y_{j+1})}$ by $[\sum_{x_{j+1}} \ee^{-V(x_{j+1})}]^2$. Under the additional assumption that $\E\{ [\sum_{|x|=1}1 ]^2\} + \E\{ [\sum_{|x|=1}\ee^{-(1+\delta)V(x)}]^2\} <\infty$ for some $\delta>0$, we have $\E\{ [\sum_{|x|=1} \ee^{-V(x)} |V(x)+b+1|^3]^2\} \le c_8\, (b+1)^6$ for some $c_8>0$ and all $b\ge 0$. Therefore,
\begin{eqnarray*}
    \Lambda_n
 &\le&c_9 (\log n)^4\, 
    \E \Big\{ \sum_{j=0}^n \sum_{|z|=j}
    {\bf 1}_{ \{ V(z_i) \in I_i, \; \forall 1\le i\le j\} } 
    \ee^{-2V(z)} (V(z)+1)^6 \Big\}
    \\
 &&+ c_{10}n^3\, 
    \E \Big\{ \sum_{j=n+1}^{2n-1} {1\over 2n-j}
    \sum_{|z|=j}
    {\bf 1}_{ \{ V(z_i) \in I_i, \; \forall 1\le i\le j\} } 
    \ee^{-2V(z)} \Big\}.
\end{eqnarray*}

\noindent Applying (\ref{many-to-one}), this leads to:
$$
\Lambda_n
\le
c_9 (\log n)^4
\sum_{j=0}^n \E \Big\{ \ee^{-S_j} (S_j+1)^6 
{\bf 1}_{ \{ \underline{S}_j \ge 0 \} } \Big\}
+
c_{10}n^3 \sum_{j=n+1}^{2n-1} {1\over 2n-j}
\E \Big\{ \ee^{-S_j}
{\bf 1}_{ \{ S_i \in I_i, \; \forall 1\le i\le j\} }\Big\}.
$$

\noindent It is easy to bound the two expectation expressions on the right-hand side. For the first expectation, we simply use Lemma \ref{l:3} to see that $\E \{ \ee^{-S_j} (S_j+1)^6 {\bf 1}_{ \{ \underline{S}_j \ge 0\} } \} \le {c_{11}\over (j+1)^{3/2}}$ for some $c_{11}>0$ and all $j\ge 0$. For the second, we recall that for $j\in [n+1, \, 2n-1]\cap \z$, $S_j \in I_j$ means $S_j \ge {3\over 2}\log n$, so that
\begin{eqnarray*}
    \E \Big\{ \ee^{-S_j}
    {\bf 1}_{ \{ S_i \in I_i, \; \forall 1\le i\le j\} }\Big\}
 &\le & \E \Big\{ \ee^{-S_j}
    {\bf 1}_{ \{ \underline{S}_{j-1} \ge 0, \;
    {3\over 2}\log n \le S_j < 3\log n\} }\Big\} + {1\over n^3}
    \\
 &\le & {1\over n^{3/2}} \,
    \P \Big\{ \underline{S}_{j-1} \ge 0, \;
    {3\over 2}\log n \le S_j < 3\log n \Big\} + {1\over n^3}
    \\
 &\le& {1\over n^{3/2}} \, c_{19} \, {9(\log n)^2 \over j^{3/2}}
    + {1\over n^3},
\end{eqnarray*}

\noindent the last inequality being a consequence of Lemma \ref{l:3}. Accordingly,
$$
\Lambda_n
\le
c_9 c_{11} (\log n)^4 \sum_{j=0}^n {1\over (j+1)^{3/2}} + c_{10} n^3 \sum_{j=n+1}^{2n-1} {1\over 2n-j} ( {9c_{19} (\log n)^2 \over n^{3/2} j^{3/2}} + {1\over n^3}) \le c_{12} (\log n)^4.
$$

\noindent Since $\E(Z_n^2) = \E(Z_n) + \Lambda_n$, and $\E(Z_n) \ge c_4$, this yields $\E(Z_n^2) \le c_{13} (\log n)^4\, [\E(Z_n)]^2$. By the Cauchy--Schwarz inequality, $\P\{ Z_n >0\} \ge {[\E(Z_n)]^2 \over \E(Z_n^2)} \ge {1\over c_{13} (\log n)^4}$; hence
$$
    \P \Big\{ \min_{|x|=2n} V(x) \le {3\over 2}\log n +2C \Big\}
    \ge
    \P\{ Z_n >0\}
    \ge
    {1\over c_{13} (\log n)^4}.
$$

\noindent We obviously can apply the same argument to study $\min_{|x|=2n-1} V(x)$, to see that for some constants $\widetilde{C}>0$ and $c_{14}>0$, and all $n\ge 2$,
\begin{equation}
    \P \Big\{ \min_{|x|=n} V(x) \le {3\over 2}\log n +\widetilde{C}\Big\}
    \ge
    {1\over c_{14} (\log n)^4}.
    \label{P()>}
\end{equation}

Let\footnote{See the last Remark in the Introduction. {F}rom here, the proof is routine, following McDiarmid~\cite{mcdiarmid}.} $\varepsilon>0$ and let $\tau_n := \inf\{ k: \; \# \{ x: \, |x|=k\} \ge n^\varepsilon\}$. For all large $n$,
\begin{eqnarray*}
 &&\P\Big\{ \tau_n <\infty, \;
    \max_{k\in [{n\over 2}, \, n]}
    \min_{|x|= k+\tau_n} V(x)
    > \max_{|y|=\tau_n} V(y) + {3\over 2}\log n +\widetilde{C}\Big\}
    \\
 &\le&\sum_{k\in [{n\over 2}, \, n]}
    \P\Big\{ \tau_n <\infty, \;
    \min_{|x|= k+\tau_n} V(x)
    > \max_{|y|=\tau_n} V(y)+ {3\over 2}\log n +\widetilde{C}\Big\}
    \\
 &\le&\sum_{k\in [{n\over 2}, \, n]}
    \Big[ \P\Big\{ \min_{|x| =k} V(x)
    > {3\over 2}\log n + \widetilde{C} \Big\}
    \Big]^{\lfloor n^\varepsilon\rfloor}  ,
\end{eqnarray*}

\noindent which, according to (\ref{P()>}), is summable in $n$. By the Borel--Cantelli lemma, a.s.\ for all sufficiently large $n$, we have either $\tau_n=\infty$, or $\max_{k\in [{n\over 2}, \, n]} \min_{|x|= k+\tau_n} V(x) \le \max_{|y|=\tau_n} V (y) +{3\over 2}\log n + \widetilde{C}$. By (\ref{LGN2}), on the system's non-extinction, a.s.\ for all large $n$, we have either $\tau_n=\infty$, or $\max_{k\in [{n\over 2}, \, n]} \min_{|x|= k+\tau_n} V(x) \le c_{15} \tau_n +{3\over 2}\log n + \widetilde{C}$.

Recall that the number of particles in each generation forms a super-critical Galton--Watson tree. In particular, on the system's non-extinction, ${\# \{ x: \, |x|= k\} \over m^k}$ converges\footnote{Here, $m := \E[\sum_{|x|=1} 1] \in (1, \, \infty)$ is the mean reproduction number in the Galton--Watson process.} a.s.\ to a positive random variable when $k\to \infty$, which implies ${\tau_n\over \log n} \to {\varepsilon \over \log m}$, a.s.\ ($n\to \infty$), and $\max_{k\in [{n\over 2}, \, n]} \min_{|x|= k+\tau_n} V(x) \ge \min_{|x|= n} V(x)$ a.s.\ for all large $n$. As a consequence, upon the system's survival, we have, a.s.\ for all large $n$,
$$
\min_{|x|= n} V(x) \le ({3\over 2} + {2\varepsilon c_{15} \over \log m})\log n + \widetilde{C}.
$$

\noindent Since $\varepsilon>0$ can be as small as possible, this yields: on the set of non-extinction,
$$
\limsup_{n\to \infty} \, {1\over \log n} \min_{|x|=n} V(x)
\le
{3\over 2}, \qquad\hbox{\rm a.s.}
$$

\noindent A fortiori, we obtain (\ref{ub}).\hfill$\Box$

\appendix
\section{Appendix on sums of i.i.d.\ random variables}

$\phantom{aob}$We list a few elementary properties of one-dimensional random walks needed in this note; they are either known results in the literature, or simple consequences of known results. Let $S_1$, $S_2-S_1$, $S_3-S_2$, $\cdots$ be i.i.d.\ real-valued random variables such that $\E(S_1)=0$ and that $0<\E(S_1^2)<\infty$. A trivial consequence of Stone's local limit theorem is that there exist constants $c_{16} >0$  and $C_0>0$ such that
\begin{equation}
    \sup_{r\in \r}
    \P\{ r\le S_n \le r+h \}
    \le
    c_{16} \, {h\over n^{1/2}},
    \qquad \forall n\ge 1, \;
    \forall h\ge C_0.
    \label{A1}
\end{equation}

\noindent We also recall (see Kozlov~\cite{kozlov}) two well-known estimates for the tail behaviour of $\underline{S}_n := \min_{1\le i\le n}S_i$: for some constant $c_{17}>0$,
\begin{eqnarray}
    \P \Big\{ \underline{S}_n \ge 0 \Big\}
 &\sim& {c_{17}\over n^{1/2}},
    \qquad n\to \infty,
    \label{A2}
    \\
    \limsup_{n\to \infty} \,
    n^{1/2} \sup_{u\ge 0} {1\over u+1} \,
    \P \Big\{ \underline{S}_n \ge -u \Big\}
 &<& \infty.
    \label{A3}
\end{eqnarray}

\medskip

\begin{lemma}
 \label{l:3}
 Let $C_0>0$ be the constant in (\ref{A1}). 
 There exists $c_{18}>0$ such that for 
 $a\ge 0$, $b\ge -a$ and
 $n\ge 1$,
 $$
     \P\Big\{ b\le S_n \le b+C_0, \;
     \underline{S}_n \ge -a\Big\}
     \le
     c_{18} \, {[(a+1)\wedge n^{1/2}] \, [(b+a+1) \wedge n^{1/2}] \over n^{3/2}} ,
 $$
 where $x\wedge y := \min \{ x, \, y\}$. In particular, there exists
 $c_{19}>0$ such that for $a\ge 0$, $b\ge -a$ and
 $n\ge 1$,
 $$
     \P\Big\{ S_n \le b, \;
     \underline{S}_n \ge -a\Big\}
     \le
     c_{19} \, {[(a+1)\wedge n^{1/2}] \, [(b+a+1)^2 \wedge n] \over n^{3/2}} .
 $$

\end{lemma}

\medskip

\noindent {\it Proof.} We only need to prove the first inequality.

There is nothing to prove if $n\le 99$; so let us assume that $n\ge 100$. We present the proof only for the case that $n$ is a multiple of $3$; say $n=3k$. A similar argument applies if $n=3k+1$ or if $n=3k+2$.

By the Markov property at time $k$, we have
\begin{eqnarray*}
 &&\P \Big\{ b\le S_{3k} \le b+C_0, \; \underline{S}_{3k} \ge -a \Big\}
    \\
 &\le&\P \Big\{ \underline{S}_k \ge -a \Big\} \sup_{x\ge -a} \P \Big\{ b-x\le S_{2k} \le b-x + C_0, \; \underline{S}_{2k} \ge -a-x \Big\}.
\end{eqnarray*}

\noindent By (\ref{A3}), $\P\{ \underline{S}_k \ge -a \} \le c_{20} {(a+1) \wedge k^{1/2}\over k^{1/2}}$. It remains to check that
$$
\sup_{x\ge -a} \P \Big\{ b-x\le S_{2k} \le b-x + C_0, \; \underline{S}_{2k} \ge -a-x \Big\} 
\le 
c_{21} \, {(b+a+1) \wedge k^{1/2}\over k}.
$$

\noindent Let $\widetilde{S}_j := S_{2k-j} -S_{2k}$. Then $\P \{ b-x\le S_{2k} \le b-x + C_0, \; \underline{S}_{2k} \ge -a-x \} \le \P \{ -b+x-C_0 \le \widetilde{S}_{2k} \le -b+x, \; \min_{1\le i\le 2k} \widetilde{S}_i \ge -a-b-C_0 \}$. By the Markov property, this leads to: for $x\ge -a$,
\begin{eqnarray*}
 &&\P \Big\{ b-x\le S_{2k} \le b-x + C_0, \; \underline{S}_{2k} \ge -a-x \Big\}
    \\
 &\le&\P \Big\{ \min_{1\le i\le k} \widetilde{S}_i \ge -a-b-C_0 \Big\} \sup_{y\in \r} \P \Big\{ -b+x-C_0-y \le \widetilde{S}_k \le -b+x-y \Big\}.
\end{eqnarray*}

\noindent The first probability expression on the right-hand side is bounded by a constant multiple of ${(b+a+1) \wedge k^{1/2} \over k^{1/2}}$ (by (\ref{A3})), whereas the second probability expression bounded by a constant multiple of $k^{-1/2}$ (by (\ref{A1})). Lemma \ref{l:3} is proved.\footnote{We mention that in the case $a=0$, Lemma \ref{l:3} is essentially Lemma 20 of Vatutin and Wachtel~\cite{vatutin-wachtel}.}\hfill$\Box$

\medskip

\begin{lemma}
 \label{l:1}
 There exists a constant
 $C>0$ such that for any $0<a\le b<\infty$,
 $$
     \liminf_{n\to \infty}
     n^{1/2} \,
     \inf_{u\in [a n^{1/2}, \, bn^{1/2}]}
     \P\Big\{ u \le S_n < u+ C
     \, \Big| \, \underline{S}_n \ge 0\Big\}
     >0.
 $$

\end{lemma}

\medskip

\noindent {\it Proof.} Follows immediately from a conditional local limit theorem (Caravenna~\cite{caravenna}): if the distribution of $S_1$ is non-lattice (i.e., not supported in any $a+b\z$, with $a\in \r$ and $b>0$), then\footnote{For any $r\in \r$, $r^+ := \max\{ r, \, 0\}$ denotes its positive part.} for any $h>0$, $\P\{ r \le S_n \le r+h \, | \, \underline{S}_n \ge 0\} = {h r^+\over n \E(S_1^2)}\exp(- {r^2\over 2n\E(S_1^2)}) + o({1\over n^{1/2}})$, $n\to \infty$, uniformly in $r\in \r$; if the distribution of $S_1$ is lattice, and is supported in $a+b\z$ with $b>0$ being the largest such value (called the ``span" in the literature), then $\P\{ S_n = a n + b \ell \, | \, \underline{S}_n \ge 0\} = {b( an + b\ell )^+ \over n\E(S_1^2)}\exp(- {(an + b\ell )^2\over 2 \E(S_1^2)}) + o({1\over n^{1/2}})$, $n\to \infty$, uniformly in $\ell \in \z$.\hfill$\Box$

\medskip

\begin{lemma}
 \label{l:4}
 Let $C>0$ be the constant in Lemma \ref{l:1}.
 For any sequence $(a_n)$ of non-negative
 numbers such that
 $\limsup_{n\to \infty} {a_n \over n^{1/2}}<\infty$, we have
 \begin{equation}
     \liminf_{n\to \infty} \, n^{3/2} \,
     \P \Big\{ \underline{S}_n \ge 0, \;
     \min_{n<j\le 2n} S_j \ge a_n, \; a_n \le S_{2n} \le a_n + 2C\Big\}
     >
     0.
     \label{A4}
 \end{equation}

\end{lemma}

\medskip

\noindent {\it Proof.} Let $c_{22}>0$ and $n_0\ge 1$ be such that $a_n \le c_{22} n^{1/2}$, $\forall n\ge n_0$. Let $p_n$ denote the probability in (\ref{A4}). Writing $\lambda_k := 2c_{22}\, n^{1/2} +kC$ for $k\ge 0$, we have, for $n\ge n_0 + \lceil ({C\over c_{22}})^2\rceil$,
\begin{eqnarray*}
    p_n
 &\ge& \sum_{k=0}^{\lfloor n^{1/2}\rfloor}
    \P \Big\{ \underline{S}_n \ge 0, \;
    \lambda_k \le S_n<\lambda_{k+1}, \;
    \\
 &&\qquad
    \min_{n<j\le 2n} (S_j-S_n) \ge a_n - \lambda_k, \;
    a_n-\lambda_k \le S_{2n}-S_n \le a_n + 2C-\lambda_{k+1}\Big\} .
\end{eqnarray*}

\noindent Note that $2C-\lambda_{k+1} = C-\lambda_k$. By independence, the probability on the right-hand side is
$$
\P \Big\{ \underline{S}_n \ge 0, \;
    \lambda_k\le S_n < \lambda_{k+1} \Big\}
    \, \P \Big\{ \underline{S}_n \ge a_n - \lambda_k, \;
    a_n-\lambda_k \le S_n \le a_n + C-\lambda_k \Big\} .
$$

\noindent The first probability expression is, by (\ref{A2}) and Lemma \ref{l:1}, greater than ${c_{23}\over n}$ (for large $n$), uniformly in $0\le k\le \lfloor n^{1/2}\rfloor$, whereas the second is, by writing $\widehat{S}_j:= S_{n-j}-S_n$, $\ge \P\{ \min_{1\le j\le n} \widehat{S}_j \ge a_n - \lambda_k +C, \; -a_n- C+\lambda_k\le \widehat{S}_n \le -a_n +\lambda_k\}$, which is $\ge \P\{ \min_{1\le j\le n} \widehat{S}_j \ge 0, \; -a_n- C+\lambda_k\le \widehat{S}_n \le -a_n +\lambda_k\}$, and thus by (\ref{A2}) and Lemma \ref{l:1} again, greater than ${c_{24}\over n}$ (for large $n$), uniformly in $0\le k\le \lfloor n^{1/2}\rfloor$. Consequently, for all sufficiently large $n$, $p_n \ge \sum_{k=0}^{\lfloor n^{1/2}\rfloor} {c_{23}\over n} \times {c_{23}\over n}$, proving the lemma.\hfill$\Box$

\bigskip
\bigskip


{\footnotesize

\baselineskip=12pt

\noindent
\begin{tabular}{lll}
& \hskip10pt Elie A\"\i d\'ekon
    & \hskip10pt Zhan Shi \\
& \hskip10pt Department of Mathematics and Computer Science
    & \hskip10pt Laboratoire de Probabilit\'es UMR
7599 \\
& \hskip10pt Technische Universiteit Eindhoven
    & \hskip10pt Universit\'e Paris VI \\
& \hskip10pt P.O.~Box 513
    & \hskip10pt 4 place Jussieu\\
& \hskip10pt 5600 MB Eindhoven
    & \hskip10pt F-75252 Paris Cedex 05\\
& \hskip10pt The Netherlands
    & \hskip10pt France \\
& \hskip10pt {\tt elie.aidekon@gmail.com}
    & \hskip10pt
    {\tt zhan.shi@upmc.fr}
\end{tabular}

}

\end{document}